\documentclass{article}
\usepackage{amsthm,amsfonts,amsbsy, amssymb,amsmath,graphicx}
\usepackage{graphics}

\newtheorem{ex}{Exercise}
\newtheorem{thm}{Theorem}
\newtheorem{dfn}{Definition}
\newtheorem{lm}{Lemma}
\newtheorem{re}{Remark}
\newtheorem{st}{Statement}
\newtheorem{prop}{Proposition}

\def\picill#1by#2(#3)
{\vbox to #2
{\hrule width #1 height 0pt depth 0pt
\vfill\epsffile{#3}}}

\makeatletter
\long\def\UR#1{\leavevmode\setbox\@tempboxa\hbox{#1}\@tempdima\fboxrule
     \advance\@tempdima \fboxsep \advance\@tempdima \dp\@tempboxa
    \hbox{\lower \@tempdima\hbox
   {\vbox{\hrule \@height \fboxrule
           \hbox{  \hskip\fboxsep
           \vbox{\vskip\fboxsep \box\@tempboxa\vskip\fboxsep}\hskip
                  \fboxsep\vrule \@width \fboxrule}%
                   }}}}
\makeatother

\makeatletter
\long\def\LR#1{\leavevmode\setbox\@tempboxa\hbox{#1}\@tempdima\fboxrule
     \advance\@tempdima \fboxsep \advance\@tempdima \dp\@tempboxa
    \hbox{\lower \@tempdima\hbox
   {\vbox{
           \hbox{  \hskip\fboxsep
           \vbox{\vskip\fboxsep \box\@tempboxa\vskip\fboxsep}\hskip
                  \fboxsep\vrule \@width \fboxrule}%
                  \hrule \@height \fboxrule}}}}
\makeatother

\makeatletter
\long\def\UL#1{\leavevmode\setbox\@tempboxa\hbox{#1}\@tempdima\fboxrule
     \advance\@tempdima \fboxsep \advance\@tempdima \dp\@tempboxa
    \hbox{\lower \@tempdima\hbox
   {\vbox{\hrule \@height \fboxrule
           \hbox{\vrule \@width \fboxrule \hskip\fboxsep
           \vbox{\vskip\fboxsep \box\@tempboxa\vskip\fboxsep}\hskip
                  \fboxsep }%
                   }}}}
\makeatother

\makeatletter
\long\def\LL#1{\leavevmode\setbox\@tempboxa\hbox{#1}\@tempdima\fboxrule
     \advance\@tempdima \fboxsep \advance\@tempdima \dp\@tempboxa
    \hbox{\lower \@tempdima\hbox
   {\vbox{
           \hbox{\vrule \@width \fboxrule \hskip\fboxsep
           \vbox{\vskip\fboxsep \box\@tempboxa\vskip\fboxsep}\hskip
                  \fboxsep }%
                  \hrule \@height \fboxrule}}}}
\makeatother

\let \ttorg \tt \def \tt{\ttorg \obeyspaces}

\date{}

 \author{Louis H. Kauffman \\
 and\\
 Vassily O. Manturov}

\title{Virtual Knots and Links}

\begin{document}

\maketitle

\abstract{This paper is an introduction to the subject of virtual
knot theory, combined with a discussion of some specific new
theorems about virtual knots. The new results are as follows: We
prove, using  a 3-dimensional topology approach that if a
connected sum of two virtual knots $K_{1}$ and $K_{2}$ is trivial,
then so are both $K_{1}$ and $K_{2}$. We establish an algorithm,
using Haken-Matveev technique,  for recognizing virtual links.
This paper may be read as both an introduction and as a research
paper.

For more about Haken-Matveev theory and its application to
classical knot theory, see \cite{Ha,Hem,Mat,HL}.

\section{Introduction}
Virtual knot theory was proposed by Louis Kauffman in 1996, see
\cite{KaV}. The combinatorial notion of virtual knot \footnote{In
the sequel, we use the generic term ``knot'' for both knots and
links, unless otherwise specified} is defined as an equivalence
class of 4-valent plane diagrams (4-regular plane graphs with
extra structure) where a new type of
crossing (called virtual) is allowed.
\bigbreak

This theory can be regarded as a ``projection'' of the knot theory
in thickened surfaces $S_{g}\times {\bf R}$ (for example, as
studied in \cite{JKS}). Regarded from this point of view, virtual
crossings appear as artifacts of the diagram projection from $S_g$
to ${\bf R}^{2}$. In such a virtual projection diagram, one does
not know the genus of the surface from which the projection was
made, and one wants to have intrinsic rules for handling the
diagrams. The rules for handling the virtual diagrams can be
motivated  (in \cite{KaV}) by the idea of a knot diagram to its
oriented Gauss code. A Gauss code for a knot is list of crossings
encountered on traversing the knot diagram, with the signs of the
crossings indicated, and whether they are over or under in the
course of the traverse. Each crossing is encountered twice in such
a traverse, and thus the Gauss code has each crossing label
appearing twice in the list. One can define Reidemeister moves on
the Gauss codes, and thus abstract the knot theory from its planar
diagrams to a theory of all oriented Gauss codes, including those
that do not have embedded realizations in the plane. Such
non-planar Gauss codes are then described by planar diagrams with
extra ``virtual'' crossings, or by knots in thickened surfaces of
higher genus. Virtual knot theory is the theory of such Gauss
codes, not necessarily realizable in the plane. When one takes
such a non-realizable code, and attempts to draw a planar diagram,
virtual crossings are needed to complete the connections in the
plane. These crossings are artifacts of the planar projection.  It
turns out that these rules describe embeddings of knots and links
in thickened surfaces, stabilized by the addition and subtraction
of empty handles (i.e. the addition and subtraction of thickened
1-handles from the surface that do not have any part of the knot
or link embedded in them) \cite{KaV2,KaV4,Ma1,Ma8,Ma10,CKS,KUP}.

Another approach to Gauss codes for knots and links is the use
of Gauss diagrams as in
\cite{GPV}). In this paper by Goussarov, Polyak and Viro, the
virtual knot theory, taken as all Gauss diagrams up to
Reidemeister moves, was used to analyze the structure of Vassiliev
invariants for classical and virtual knots. In both \cite{KaV} and
\cite{GPV} it is proved that if two classical knots are equivalent
in the virtual category \cite{KUP}, then they are equivalent in
the classical category. Thus classical knot theory is properly
embedded in virtual knot theory. \bigbreak

To date, many invariants of classical knots have been generalized
for the virtual case, see \cite{GPV,KaV,KR,Ma1,Ma2,Ma8,Ma10,Saw,
SW}. In many cases, a classical invariant extends to an invariant
of virtual knots. In some cases one has an invariant of virtuals
that is an extension of ideas from classical knot theory that
vanishes or is otherwise trivial for classical knots. See
\cite{Saw},\cite{SW}, \cite{KR}, \cite{Ma2,Ma3,Ma5}. Such
invariants are valuable for the study of virtual knots, since they
promise the possibility of distinguishing classical from virtual
knots in key cases. On the other hand, some invariants evaluated
on classical knots coincide with well known classical knot
invariants (see \cite{KaV, KaV2, KaV4, Ma3} on generalizations of
the Jones polynomial, fundamental group, quandle and quantum link
invariants). These invariants exhibit interesting phenomena on
virtual knots and links: for instance, there exists a virtual knot
$K$ with ``fundamental group'' isomorphic to ${\bf Z}$ and Jones
polynomial not equal to $1$. This phenomenon immediately implies
that the knot $K$ is not classical, and underlines the difficulty
of extracting the Jones polynomial from the fundamental group in
the classical case. \bigbreak

We know in principle that the fundamental group, plus peripheral
information, determines the knot itself in the classical case. It
is not known how to extract the Jones polynomial from this
algebraic information. The formally defined fundamental group of a
virtual knot can be interpreted as the fundamental group of the
complement of the virtual knot in the one-point suspension of a
thickened surface where this knot is presented. \bigbreak

Another phenomenon that does not appear in the classical case are
long knots \cite{Ma11}: if we break a virtual knot diagram at two
different points and take them to the infinity, we may obtain two
different long knots. \bigbreak

In the present paper, we are going to discuss both algebraic and
geometric properties of virtual knots. We recommend as survey
papers for virtual knot theory the following
\cite{KaV2,Ma1,KUP,KaV4,KaV,FKM,FJK}. \bigbreak

The paper is organized as follows. First, we give definitions and
recall some known results. In the second section, we prove, using
a 3-dimensional topology approach that if a connected sum of two
virtual knots $K_{1}$ and $K_{2}$ is trivial, then so are both
$K_{1}$ and $K_{2}$. Here we say ``a connected sum'' because the
connected sum is generally not well defined. In the third section
we introduce Haken-Matveev theory of normal surfaces in order to
establish an algorithm for recognizing virtual knots. The fourth
section is devoted to self-linking coefficient that generalizes
the writhe of a classical knot. The fifth section gives a quick
survey of the relationship between virtual knots, welded knots and
embeddings of tori in four-dimensional space.

We do not touch many different subjects on virtual knots. For
instance, we do not describe virtual links and their invariatns.
For more details see \cite{Kam,Ver,FRR,MBraids,MBook,KL}.

 \bigbreak

\subsection{Basic definitions}

Let us start with the definitions and introduce the notation.

\begin{dfn}
A {\em virtual link diagram} \index{Virtual link diagram} is a
planar graph of valency four endowed with the following structure:
each vertex either has an overcrossing and undercrossing or is
marked by a virtual crossing, (such a crossing is shown in Fig.
\ref{vcross}).

\begin{figure}
\begin{center}
\unitlength 0.7mm
\begin{picture}(45,38)
\put(25,25){\circle{3}} \put(10,10){\line(1,1){30}}
\put(40,10){\line(-1,1){30}}
\end{picture}
\end{center}
\vspace{-1cm} \caption{Virtual crossing} \label{vcross}
\end{figure}
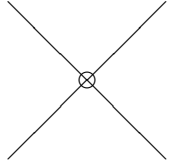

\end{dfn}

All crossings except virtual ones are said to be {\em classical}.

Two diagrams of virtual links (or, simply, {\em virtual diagrams})
are said to be {\em equivalent} if there exists a sequence of {\em
generalized Reidemeister moves}\index{Reidemeister moves
generalized}, transforming one diagram to the other one.

As in the classical case, all moves are thought to be performed
inside a small domain; outside this domain the diagram does not
change.

\begin{dfn}
Here we give the list of {\em generalized Reidemeister moves}:
\begin{enumerate}

\item Classical Reidemeister moves related to classical vertices.

\item Virtual versions $\Omega'_{1},\Omega'_{2},\Omega'_{3}$ of
Reidemeister moves, see Fig. \ref{vver}.

\begin{figure}
\begin{center}

\unitlength 0.5mm \linethickness{0.4pt}
\begin{picture}(136.67,70)
\bezier{40}(5,50)(10.33,50.33)(14.33,51.67)
\bezier{32}(14.33,51.67)(18.33,53)(22.33,51.67)
\bezier{52}(22.33,51.67)(28,50.33)(35,50)
\bezier{40}(45,50)(50.33,50.33)(54.33,51.67)
\bezier{52}(62.33,51.67)(68,50.33)(75,50)
\bezier{24}(54.33,51.67)(57,52.67)(58.67,55)
\bezier{24}(58.67,55)(61,56.33)(61,60)
\bezier{32}(61,60)(61.33,64)(57.67,63.67)
\bezier{28}(57.67,63.67)(54.33,63)(55,59.67)
\bezier{28}(55,59.67)(55,55.67)(58,54)
\bezier{20}(58,54)(59.33,52.33)(62.33,51.67)
\put(58,54){\circle{2}}
\bezier{12}(85,70)(84.67,68.67)(85.67,67.33)
\bezier{64}(85.67,67.33)(86.67,60.33)(85.33,52)
\bezier{28}(85.33,52)(84.33,48.67)(85,45)
\bezier{20}(100,45)(99,47.67)(99.67,50)
\bezier{56}(99.67,50)(99,57.33)(99.67,64.33)
\bezier{24}(99.67,64.33)(100,70)(100,70)
\bezier{12}(120,70)(119.67,68.67)(120.67,67.33)
\bezier{28}(120.33,52)(119.33,48.67)(120,45)
\bezier{20}(135,45)(134,47.67)(134.67,50)
\bezier{24}(134.67,64.33)(135,70)(135,70)
\bezier{44}(121,67.33)(128,64)(129.67,61.67)
\bezier{44}(129.67,61.67)(131.33,56)(126.67,55.33)
\bezier{32}(126.67,55.33)(122,55)(120.33,52)
\bezier{16}(134.67,50.33)(136.67,51.67)(135,52.67)
\bezier{36}(135,52.67)(128,53.33)(126.67,55.33)
\bezier{20}(126.67,55.33)(124,57)(124.33,58.33)
\bezier{44}(124.33,58.33)(122,62.67)(127.67,64.33)
\bezier{20}(134.33,64)(133.67,62)(131.67,63.67)
\bezier{32}(131.67,63.67)(129.33,67.33)(127,64)
\put(126.67,64.33){\circle{2}} \put(126.67,55){\circle{2}}
\bezier{200}(65,40)(53.33,22.67)(30,5)
\bezier{200}(30,40)(49.33,25.33)(65,5) \put(50,22.33){\circle{2}}
\bezier{200}(120,40)(108.33,22.67)(85,5)
\bezier{200}(85,40)(104.33,25.33)(120,5)
\put(105,22.33){\circle{2}}
\bezier{264}(70,22.67)(47.33,46.67)(25,22.67)
\bezier{276}(80,22.67)(102.33,-3.33)(125,22.67)
\put(114,12.67){\circle{2}} \put(93.67,11.67){\circle{2}}
\put(59,31.67){\circle{2}} \put(38.67,32.67){\circle{2}}
\put(75.33,22.67){\makebox(0,0)[cc]{$\Longleftrightarrow$}}
\put(39.67,52.33){\makebox(0,0)[cc]{$\Longleftrightarrow$}}
\put(110.67,59){\makebox(0,0)[cc]{$\Longleftrightarrow$}}
\end{picture}

\end{center}
\vspace{-0.7cm}

\caption{Moves $\Omega'_{1},\Omega'_{2},\Omega'_{3}$} \label{vver}
\end{figure}
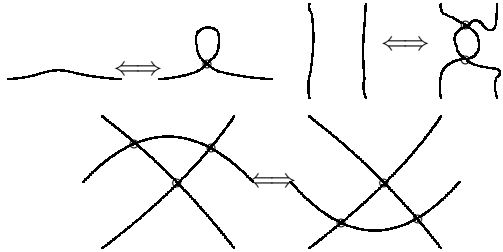

\item The ``semivirtual'' version of the third Reidemeister move,
see Fig.~\ref{semvir},

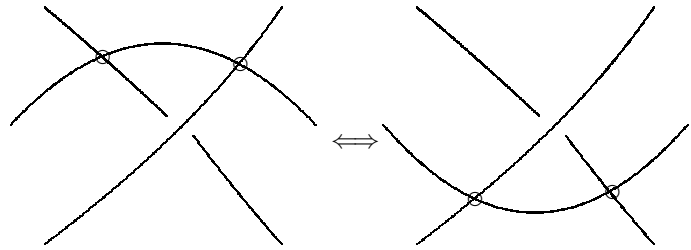
\begin{figure}
\unitlength 0.9mm
\begin{picture}(136,40)

\bezier{200}(65,40)(53.33,22.67)(30,5)
\bezier{200}(30,40)(39,32.5)(48,24)
\bezier{100}(52,21)(58,13)(65,5)
\bezier{200}(120,40)(108.33,22.67)(85,5)

\bezier{200}(85,40)(94,32.5)(103,24)
\bezier{100}(107,21)(113,13)(120,5)

\bezier{264}(70,22.67)(47.33,46.67)(25,22.67)
\bezier{276}(80,22.67)(102.33,-3.33)(125,22.67)
\put(114,12.67){\circle{2}} \put(93.67,11.67){\circle{2}}
\put(59,31.67){\circle{2}} \put(38.67,32.67){\circle{2}}

\put(76,20){\makebox(0,0)[cc]{$\Longleftrightarrow$}}
\end{picture}

\vspace{-0.3cm}

\caption{The semivirtual move $\Omega''_{3}$} \label{semvir}
\end{figure}

\end{enumerate}

\end{dfn}

\begin{re}
The two moves shown in Fig. \ref{forbid} are {\em
forbidden},\index{Forbidden move} i.e., they are not in the list
of generalized moves and cannot be expressed via these moves.

\begin{figure}
\centering\includegraphics[width=320pt, height=75pt]{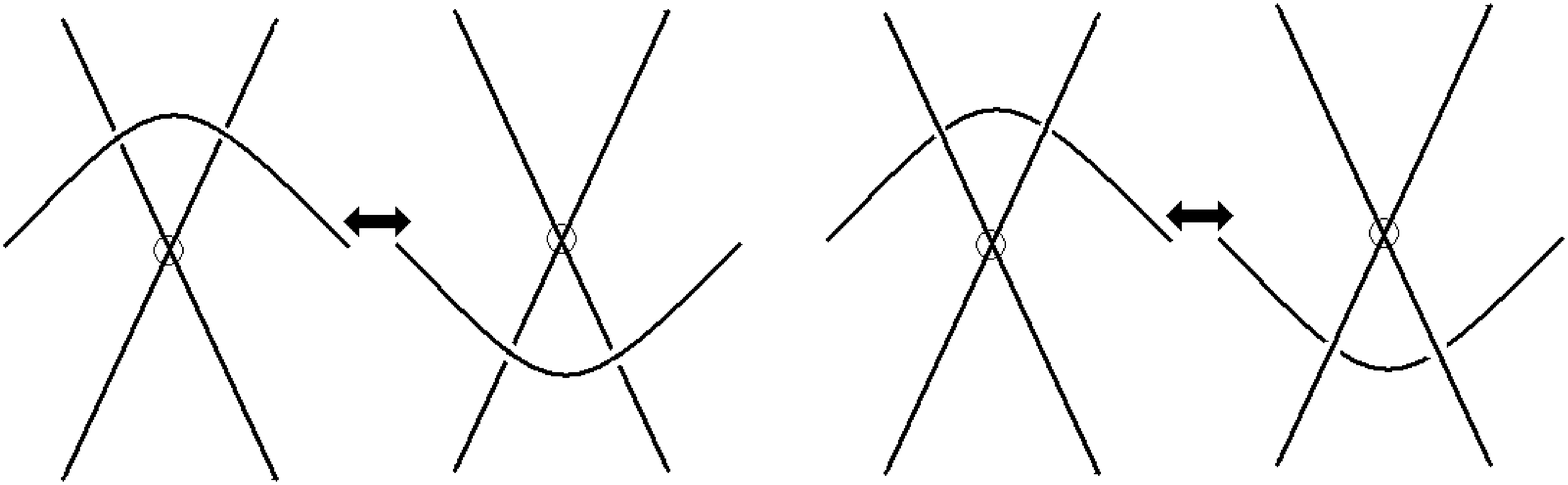}
\caption{The forbidden moves} \label{forbid}
\end{figure}

\end{re}

\begin{dfn}
A {\em virtual link} \index{Link virtual} is an equivalence class
of virtual diagrams modulo generalized Reidemeister moves.
\end{dfn}

One can easily calculate the number of components of a virtual
link. A {\em virtual knot}\index{Knot virtual} is a one--component
virtual link.

\begin{ex}
Show that any virtual link having a diagram without classical
crossings is equivalent to a classical unlink.
\end{ex}

\begin{re}
Formally, generalized Reidemeister moves appear to give a new
equivalence relation for classical links: there seem to exist two
isotopy relations for classical links, the classical one that we
are used to work with and the virtual one. In \cite{GPV}, it is
shown that, for classical knots and links, virtual and classical
equivalences are the same.
\end{re}

\begin{re}
Actually, the forbidden move is a very strong one. Each virtual
knot can be transformed to any other one by using all generalized
Reidemeister moves and the forbidden move.

This was first proved in \cite{GPV} (see also \cite{Nel,Kan}).
Therefore, any two virtual knots can be transformed to each other
by a sequence of generalised Reidemeister moves and the forbidden
moves.

If we allow only the forbidden move shown in the left part of Fig.
\ref{forbid}, we obtain what are called {\em welded knots},
developed by Shin Satoh, \cite{Satoh}. Some initial information on
this theory can be found in \cite{Kam}, see also \cite{FRR}. A
quick survey of this theory is given in section 5 of the present
paper.
\end{re}

\begin{dfn}
By a mirror image \index{Mirror image} of a virtual link diagram
we mean a diagram obtained from the initial one by switching all
types of {\bf classical} crossings (all virtual crossings stay on
the same positions).
\end{dfn}

\subsubsection{Projections from thikcened surfaces}

The choice of generalized Reidemeister moves is very natural. It
is the complete list of moves that occur when considering a
generic projection of $S_{g} \times I$ to ${\bf R}^{2} \times I$
(or equivalently ${\bf R}^{3}$), i.e. thickened surface is
projected generically to thickened plane. The virtual crossings
appear in this projection as artifacts of different sheets of the
surface being projected to the single sheet below. Actual
crossings in the thickened surface are rendered as classical
crossings in the thickened plane. The other moves, namely, the
semivirtual move and purely virtual moves, are shown in
Fig.\;\ref{haba} together with the corresponding moves in
thickened surfaces.

Note that plane diagrams modulo classical Reidemeister moves lead
to the classical knot theory as well as spherical diagrams (on
$S^{2}$) modulo Reidemeister moves; in the sequel, we shall use
spherical diagrams rather than planar ones.

%\end{figure}

\begin{figure}[t]
\centering\includegraphics[width=350pt]{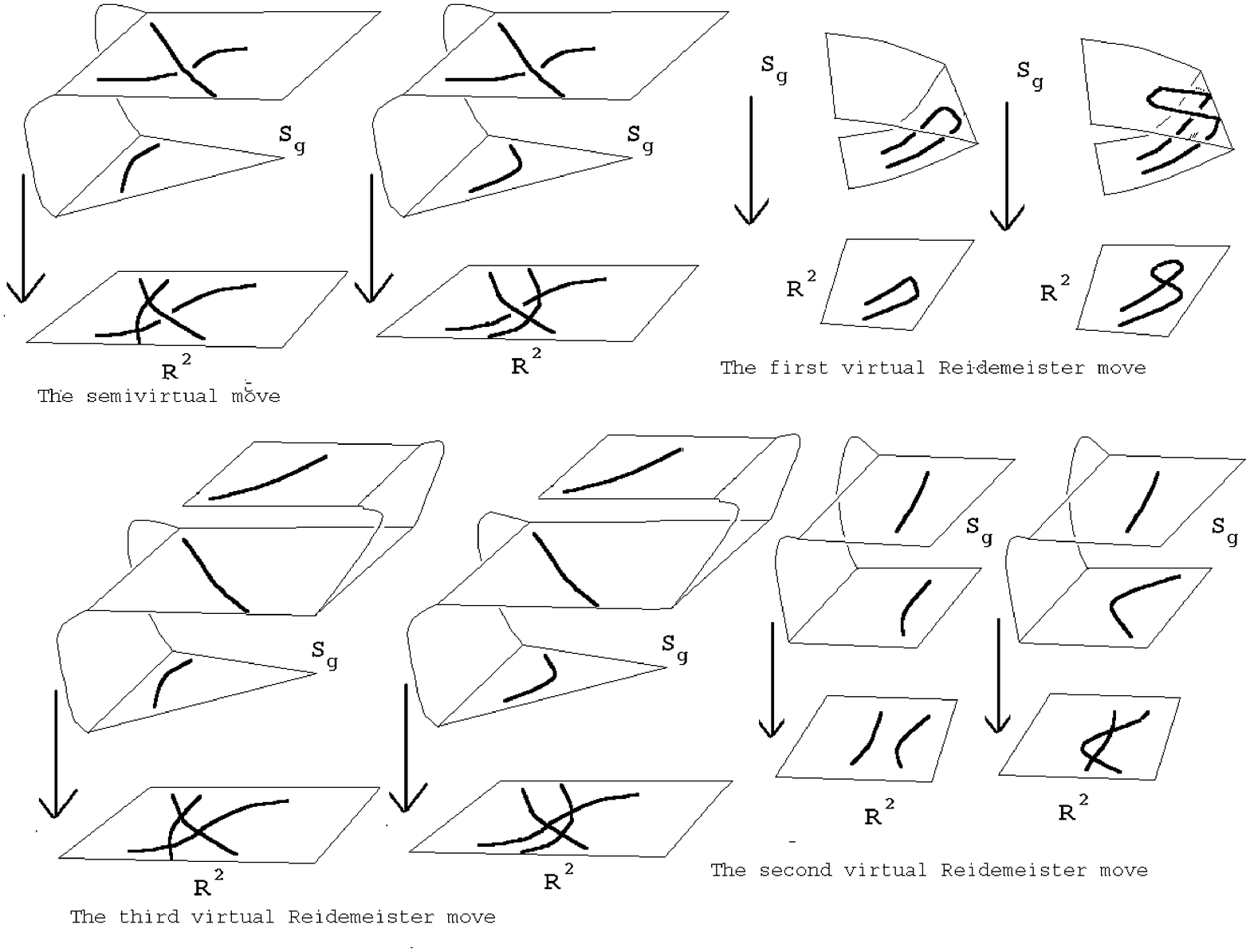}
\caption{Generalized Reidemeister moves and thickened surfaces}
\label{haba}
\end{figure}

   In fact, there exists a topological intepretation
for  virtual knot theory in terms of embeddings of links in
thickened surfaces \cite{KaV,KaV2}. Regard each virtual crossing
as a shorthand for a detour of one of the arcs in the crossing
through a $1$--handle that has been attached to the $2$--sphere of
the original diagram. The two choices for the $1$--handle detour
are homeomorphic to each other (as abstract manifolds with
boundary). By interpreting each virtual crossing in such a way, we
obtain an embedding of a collection of circles into a thickened
surface $S_{g}\times {\bf R}$, where $g$ is the number of virtual
crossings in the original diagram $L$ and $S_{g}$ is the
orientable $2$--manifold homeomorphic to the sphere with $g$
handles. Thus, to each virtual diagram $L$ we obtain an embedding
$s(L)\to S_{g(L)}\times {\bf R}$, where $g(L)$ is the number of
virtual crossings of $L$ and $s(L)$ is a disjoint union of
circles. We say that two such surface embeddings are {\em stably
equivalent} if one can be obtained from the other by isotopy in
the thickened surface, homeomorphisms of surfaces, and the
addition of substraction or handles not incident to images of
curves.

\begin{thm}
Two virtual link diagrams generate equivalent (isotopic) virtual
links if and only if their corresponding surface embeddings are
stably equivalent.
\end{thm}

This result was sketched in \cite{KaV}. The complete proof appears
in \cite{KaV3}.

A hint to this proof is demonstrated in Fig. \ref{haba}.

Here we wish to emphasize the following important circumstance.

\begin{dfn}
A virtual link diagram is {\em minimal}\index{Minimal virtual link
diagram} if no handles can be removed after a sequence of
Reidemeister moves.
\end{dfn}

An important Theorem by Kuperberg \cite{KUP} says the following.

\begin{thm}
For a virtual knot diagram $K$ there exists a unique minimal
surface in which an $I$--neighbourhood of an equivalent diagram
embeds and the embedding type of the surface is unique up to
isotopy of the image in the thickened surface.\label{Kpbg}
\end{thm}

\begin{dfn}
The {\em genus}, $g(K)$, of a virtual knot or link is the genus of
the unique minimal surface described in theorem \ref{Kpbg}.
\end{dfn}

\subsubsection{Gauss diagram approach}
\begin{dfn}
A {\em Gauss diagram}\index{Gauss diagram} of a (virtual) knot
diagram $K$ is an oriented circle (with a fixed point) where
pre-images of overcrossing and undercrossing of each {\bf
classical} crossing are connected by a chord. Pre-images of each
crossing are connected by an arrow, directed from the pre-image of
the overcrossing to the pre-image of the undercrossing. The sign
of each arrow equals the local writhe number of the vertex
(defined as in the classical case). Note that arrows (chords)
correspond only to classical crossings.
\end{dfn}

\begin{re}
For classical knots this definition coincides with the standard
one.
\end{re}

Given a Gauss diagram with labelled arrows, if this diagram is
realizable then it (uniquely) represents some classical knot
diagram. Otherwise one cannot get any classical knot diagram.

Herewith, the four--valent graph represented by this Gauss diagram
and not embeddable in ${\bf R}^{2}$ can be {\em immersed} to ${\bf
R}^{2}$. Certainly, we shall consider only ``good'' immersions
without triple points and tangencies.

Having such an immersion, let us associate virtual crossings with
intersections of edge images, and classical crossings at images of
crossings, see Fig. \ref{vcclcr}.

\begin{figure}
\begin{center}
\unitlength 1mm \linethickness{0.4pt}
\begin{picture}(97.67,39)
\multiput(19.67,35)(0.99,-0.10){3}{\line(1,0){0.99}}
\multiput(22.64,34.70)(0.36,-0.11){8}{\line(1,0){0.36}}
\multiput(25.49,33.82)(0.22,-0.12){12}{\line(1,0){0.22}}
\multiput(28.12,32.39)(0.13,-0.11){17}{\line(1,0){0.13}}
\multiput(30.40,30.47)(0.12,-0.15){16}{\line(0,-1){0.15}}
\multiput(32.27,28.14)(0.11,-0.22){12}{\line(0,-1){0.22}}
\multiput(33.63,25.48)(0.12,-0.41){7}{\line(0,-1){0.41}}
\multiput(34.44,22.60)(0.11,-1.49){2}{\line(0,-1){1.49}}
\multiput(34.66,19.63)(-0.09,-0.74){4}{\line(0,-1){0.74}}
\multiput(34.29,16.66)(-0.12,-0.35){8}{\line(0,-1){0.35}}
\multiput(33.34,13.83)(-0.11,-0.20){13}{\line(0,-1){0.20}}
\multiput(31.85,11.24)(-0.12,-0.13){17}{\line(0,-1){0.13}}
\multiput(29.87,9)(-0.15,-0.11){16}{\line(-1,0){0.15}}
\multiput(27.49,7.20)(-0.24,-0.12){11}{\line(-1,0){0.24}}
\multiput(24.80,5.90)(-0.41,-0.11){7}{\line(-1,0){0.41}}
\multiput(21.90,5.17)(-1.49,-0.07){2}{\line(-1,0){1.49}}
\multiput(18.92,5.02)(-0.74,0.11){4}{\line(-1,0){0.74}}
\multiput(15.97,5.46)(-0.31,0.11){9}{\line(-1,0){0.31}}
\multiput(13.16,6.49)(-0.20,0.12){13}{\line(-1,0){0.20}}
\multiput(10.61,8.04)(-0.13,0.12){17}{\line(-1,0){0.13}}
\multiput(8.42,10.07)(-0.12,0.16){15}{\line(0,1){0.16}}
\multiput(6.68,12.50)(-0.11,0.25){11}{\line(0,1){0.25}}
\multiput(5.45,15.22)(-0.11,0.49){6}{\line(0,1){0.49}}
\put(4.78,18.13){\line(0,1){2.99}}
\multiput(4.71,21.12)(0.10,0.59){5}{\line(0,1){0.59}}
\multiput(5.23,24.06)(0.11,0.28){10}{\line(0,1){0.28}}
\multiput(6.32,26.84)(0.12,0.18){14}{\line(0,1){0.18}}
\multiput(7.94,29.35)(0.12,0.12){18}{\line(0,1){0.12}}
\multiput(10.02,31.49)(0.16,0.11){15}{\line(1,0){0.16}}
\multiput(12.49,33.17)(0.28,0.12){10}{\line(1,0){0.28}}
\multiput(15.25,34.33)(0.74,0.11){6}{\line(1,0){0.74}}
\put(4.67,20){\vector(1,0){30}}
\put(26,33.67){\vector(0,-1){27.33}}
\put(13.33,6.33){\vector(0,1){27.33}}
\put(10.33,31.67){\vector(-1,-1){1}}
\put(30.67,30.33){\circle*{1.33}}
\bezier{250}(75,20)(50.33,38.67)(75,38.67)
\bezier{100}(74.33,38.67)(86.33,39)(86.67,35.33)
\bezier{200}(64.33,34)(96.33,36)(93.67,29.67)
\put(94,30.67){\vector(-1,3){0.2}}
\bezier{100}(93.33,13.33)(97.67,22)(94,30.67)
\bezier{100}(86.67,32.67)(86,25)(75,19.67)
\bezier{100}(75,20)(64.33,14.33)(75,8.33)
\bezier{200}(72.67,6.33)(89.67,8)(75,20)
\bezier{100}(80.33,7)(89,4.33)(92.67,12.33)
\bezier{250}(61.33,33.67)(47.33,10)(72.33,6.33)
\put(93,13){\circle*{1.33}} \put(28,35.33){\makebox(0,0)[cc]{$1$}}
\put(13,37){\makebox(0,0)[cc]{$2$}}
\put(2,20){\makebox(0,0)[cc]{$3$}}
\put(9,22.33){\makebox(0,0)[cc]{$-$}}
\put(15.67,29){\makebox(0,0)[cc]{$+$}}
\put(23.67,29){\makebox(0,0)[cc]{$+$}}
\put(89,36.33){\makebox(0,0)[cc]{$1$}}
\put(83,31.67){\makebox(0,0)[cc]{$+$}} \put(75,20){\circle{2}}
\put(59,36.33){\makebox(0,0)[cc]{$2$}}
\put(67,31.33){\makebox(0,0)[cc]{$+$}}
\put(84,8.67){\makebox(0,0)[cc]{$3$}}
\put(71,8.33){\makebox(0,0)[cc]{$-$}}
\put(46.33,20.67){\makebox(0,0)[cc]{$\Longrightarrow$}}

\end{picture}
\end{center}
\vspace{-0.6cm} \caption{A virtual knot and its Gauss diagram}
\label{vcclcr}
\end{figure}
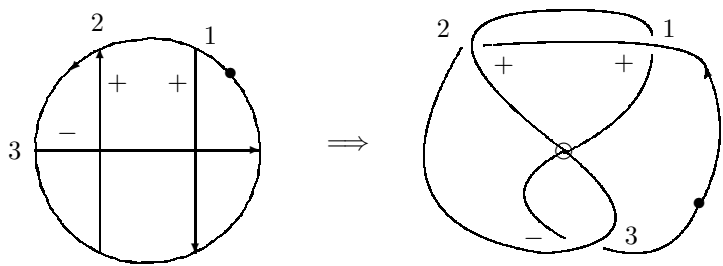

Thus, by a given Gauss diagram we have constructed (not uniquely)
a virtual knot diagram.

\begin{thm}[\cite{GPV}]
The virtual knot isotopy class is uniquely defined by this Gauss
diagram.
\end{thm}

\begin{re}
The equivalent result for Gauss codes is discussed in \cite{KaV}.
\end{re}

\begin{ex}
Prove this fact.
\end{ex}

\begin{ex}
Show that purely virtual moves and the semivirtual move are
exactly those moves that do not change the Gauss diagram at all.
\end{ex}

\section{Underlying genus of virtual knots}

The connected sum for two classical knots is defined by breaking
these diagrams at one point each and attaching one diagram to the
other one respecting the orientation. This is a well-defined
operation. Analogously, one can define a connected sum for virtual
knots by means of their diagrams. It is well known that this
operation {\bf is not} well defined, see, e.g., \cite{MBook}: the
definition of a connected sum depends on initial diagrams and the
choice of break point.

Figure \ref{Kis} illustrates a non-trivial connected sum of
trivial virtual knots. This example is due to Kishino, \cite{KS}.

\begin{figure}
\centering\includegraphics[width=300pt]{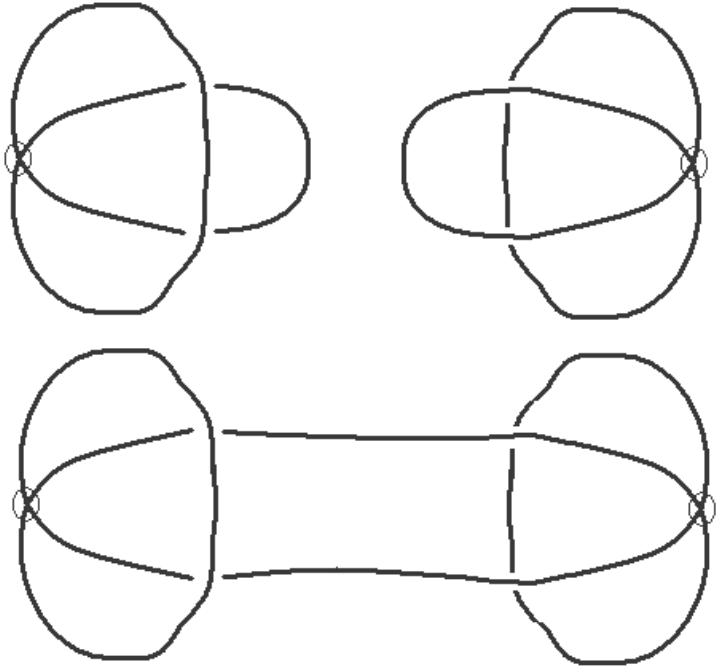} \caption{A
non-trivial connected sum of trivial virtual knots} \label{Kis}
\end{figure}

In what follows, we shall prove results for a connected sum
$K_{1}{\#}K_{2}$, with the intent that our statement holds {\em
for any} connected sum: here the notation $K_{1}{\#}K_{2}$ will be
used for an arbitrary connected sum, unless otherwise specified.

The aim of this section is to prove the following.

\begin{thm}
Let $K_{1}, K_{2}$ be two virtual knots. Then $g(K_{1}\# K_{2})\ge
g(K_{1})+g(K_{2})-1$, whereas if $g(K_{1})=0$ or $g(K_{2})=0$ then
$g(K_{1}\# K_{2})\ge g(K_{1})+g(K_2)$. Here $g(K)$ is the genus
defined at the end of the last section, the minimal genus in which
$K$ can be represented.\label{tht}
\end{thm}

From this theorem one can deduce the following

\begin{thm}
If $K_{1}$ and $K_{2}$ are virtual knots such that some connected
sum $K_{1}\# K_{2}$ is trivial, then both $K_{1}$ and $K_{2}$ are
trivial knots.\label{th1mn}
\end{thm}

Indeed, if at least one of $K_{1},K_{2}$ has positive underlying
genus, then the underlying genus of their connected sum is also
positive (by Theorem \ref{tht}). In the case when both knots have
underlying genus zero, one should mention the following.

\begin{lm}
If $K_{1}$ and $K_{2}$ are two classical knots then any connected
sum $K_{1}{\#}K_{2}$ of genus zero is equivalent to the
(well-defined) classical connected sum
$K_{1}{\#}K_{2}$.\label{ord}
\end{lm}

This lemma follows directly from the proof of Theorem \ref{tht}.

The remaining part of Theorem \ref{th1mn} now follows from the
non-triviality of connected sum in the classical case. See
\cite{CF}.

\subsection{Two types of connected sums}

Given two virtual knots $K_{1}$ and $K_{2}$ represented by knots
in thickened surfaces, there are two natural possibilities to
represent their connected sum as a knot in a thickened surface.

\begin{re}
In what follows, we shall use the same letter for an abstract
virtual knot $K$ and the knot lying in a thickened surface
$M\times I$ and representing the knot $K$ (abusing notation).
\end{re}

We shall deal with 3-manifolds (possibly, with boundary) and
2-surfaces in these manifolds. A compact surface $F$ in a manifold
$M$ is called {\em proper} if $F\cap
\partial M=\partial F$. In the sequel, all surfaces are assumed to be proper.

Also, we assume that all surfaces considered in the manifold $M$
intersect the knot transversally.

The first method for making the connected sum goes as follows. We
take thickened surfaces $(M_{1}\times I)\supset K_{1}$ and
$(M_{2}\times I)\supset K_{2}$ and cut two vertical full cylinders
$D_{i}\times I$, where $D_{i}\subset M_{i}$, such that
$(D_{i}\times I)\cap K_{i}$ is homeomorphic to an interval. Then
we paste the obtained manifolds together (by identifying $\partial
D_{1}\times I$ and $\partial D_{2}\times I$ with respect to the
orientation of manifolds, the direction of the interval $I$ and
the orientation of knot at the two gluing points) and obtain
$(M\times I)=((M_{1}\# M_{2})\times I)$ with the knot $K_{1}\#
K_{2}$ inside.

Clearly, $g(M)=g(M_{1})+g(M_{2})$.

Another way to construct the connected sum works only in some
special cases. Suppose $K_{1}$ and $K_{2}$ lie in $M_{1}\times I,
M_{2}\times I$, where both $g(M_{1})$ and $g(M_{2})$ are greater
than zero, and there exist two non-trivial (non zero-homotopic)
curves $\gamma_{1}\subset M_{1}$ and $\gamma_{2}\subset M_{2}$
such that $(\gamma_{i}\times I)\cap K_{i}$ consists of precisely
one point (note that in this case such a curve can not divide the
$2$-manifold into two parts). Then we cut the thickened surfaces
$(M_{i}\times I)$ along $\gamma_{i}\times I$ and paste them
together. Thus we obtain a manifold $M\times I$, where
$g(M)=g(M_{1})+ g(M_{2})-1$ with some connected sum $K_{1}\#
K_{2}$ lying in $M$.

It turns out that these two connected sums are the only essential
forms of connected summation for virtual knots, from which we
deduce Theorem \ref{th1mn}.

\label{scty}

\subsection{The proof plan of Theorem \ref{th1mn}}

Consider two virtual knots $K_{1}$ and $K_{2}$ and their connected
sum $K_{1}\# K_{2}$. Let us realize this connected sum by curves
in thickened surfaces by using the first method. Denote the
corresponding surfaces by $M_{1},M_{2},M_{1}\# M_{2}$, and denote
the corresponding knots by $K_{1},K_{2},K_{1}\#K_{2}$ (abusing
notation). Now, we are going to transform $(M_{1}\#M_{2})\times I$
and knots inside it.

To simplify the notation, let us use the same letters for closed
surfaces and surfaces with boundary. We will write
$M=M_{1}{\#}M_{2}$ and $M=M_{1}\cup M_{2}$. We prefer the second
notation to emphasize that both $M_{1}$ and $M_{2}$ are parts of
$M$, while the first notation will be used when $M_{1}$ and
$M_{2}$ are treated as separate manifolds.

We will perform a destabilization process on the knot obtained by
the connected sum operation described in \ref{scty}.

We are going to check that the following conditions hold during
the destablization process (all notation stays the same during the
transformation):

\begin{enumerate}

\item The ambient manifold $M$ can be divided into two parts
$M_{1}\#M_{2}$ such that $M_{i}\times I$ represents the knot
$K_{i}$ (i.e., if we close this manifold, we obtain a surface
realization of $K_{i}$).

\item The intersection $M=M_{1}\cap M_{2}$ consists of one or two
components; so $(M_{1}\times I)\cap (M_{2}\times I)$ consists of
one or two annuli.

\item The knot $K_{1}\#K_{2}$ intersects the manifold $(M_{1}\cap
M_{2})$ precisely at two points; in the case when $M_{1}\cap
M_{2}$ is not connected, these intersection points lie in
different connected components.

\item The process stops when $g(M_{1}\# M_{2})$ is the minimal
genus of the knot $K_{1}\# K_{2}$.
\end{enumerate}

If we organize the process as described above, we prove Theorem
\ref{th1mn}. Indeed, at each moment of the process we have $K_{1}$
and $K_{2}$ represented by knots in thickened surfaces of genera
$g_{1}$ and $g_{2}$. The knot $K_{1}\#K_{2}$ lies in the surface
of genus $g_{1}+ g_{2}$ if we deal with the connected sum of the
first type and in the surface of genus $g_{1}+g_{2}-1$ if we deal
with the connected sum of the second type. So, the same holds when
the process stops, thus we have $g(K_{1}{\#}K_{2})=g_{1}+g_{2}$ or
$g(K_{1}{\#}K_{2})=g_{1}+g_{2}-1$, where the last case is possible
only if we have the connected sum of the second type (hence, both
$g_{1}$ and $g_{2}$ are greater than zero). Taking into account
that $g_{i}$ is the genus of a surface (not necessarily minimal)
representing $K_{i}$, we obtain the statement of the theorem.

\subsection{The process}

In the present subsection we describe how this process works.

Suppose we have the connected sum of type $i$ ($i=1$ or $2$) of
the knots $K_{1}$ and $K_{2}$. The main statement is the
following.

\begin{st}
If there is a possibility to decrease the genus of $M_{1}{\#}
M_{2}$, then one of the following holds:

\begin{enumerate}

\item We can perform a destabilization in $M_{i}$ without changing
$M_{3-i}$ and the connected sum type (thus, we decrease the genus
of one of connected summands $M_{i}$ by one, as well as that of
$M_{1}{\#}M_{2}$).

\item If we have the first type connected sum, then there is a
possibility to transform it to the connected sum of the second
type, decreasing the genus of $M_{1} {\#} M_{2}$ by one without
changing the genera of $M_{1}$ and $M_{2}$.

\item If we have the second type connected sum, then there is a
possibility to transform it to the connected sum of the first
type, decreasing each of $g(M_{1}), g(M_{2}), g(M_{1}{\#}M_{2})$
precisely by one.

\end{enumerate}\label{sttm}
\end{st}

Together with all points described above, this statement completes
the proof of Theorem \ref{th1mn}.

\begin{proof}[Proof of Statement \ref{sttm}]
First, consider the case of the first type connected sum. We have
$M=M_{1}{\#}M_{2}$.  Denote $M_{1}\cap M_{2}$ by $D$.
 Suppose we are able to destabilize the
pair $((M_{1}{\#}M_{2})\times I, K_{1}{\#}K_{2}).$ Then there is a
vertical annulus $C$ in $M_{1}{\#} M_{2}$ which does not intersect
the knot $K_{1}{\#}K_{2}$. If there is such an annulus which does
not intersect $D$, then we can destabilize one of the summands
along $C$; this is the first case of Statement \ref{sttm}.

Suppose there is no such annulus $C$. Any $C$ we consider will
intersect $D$. Without loss of generality assume that the
intersection between each such $C$ and $D$ is transverse. Let $n$
be the minimal number of connected components of the intersection
$C\cap D$.

Since $C$ and $D$ are manifolds with boundary (vertical annuli),
their (generic) intersection  may consist of:

\begin{enumerate}

\item simple circles;

\item trivial arcs;

\item horizontal circles;

\item vertical arcs;

\end{enumerate}

Here the circle is {\em trivial} if it represents the trivial
element in the fundamental group of the annulus $C$, otherwise the
circle is called {\em horizontal}. The arc is called {\em trivial}
if it connects points from the same boundary component of the
annulus; otherwise it is called {\em vertical}.

If there is a trivial circle, then we can consider an innermost
circle $\gamma$ (with respect to $C$). This circle contains no
intersection points with $D$ inside. Because this disk together
with a disk from $D$ bounds a 3-ball, (see Fig. \ref{rm1}), we can
slightly change the annulus $C$ in such a way that the total
intersection between $C$ and $D$ decreases, and $C$ remains an
annulus with non-contractible core. The same situation happens
when we have a trivial arc, see Fig.~\ref{rm2}.

\begin{figure}
\centering\includegraphics[width=200pt]{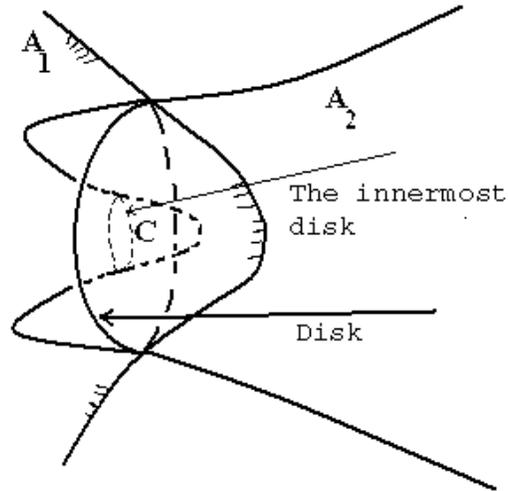}
\caption{Non-horizontal disks in $A_{1}$} \label{rm1}
\end{figure}

\begin{figure}
\centering\includegraphics[width=300pt]{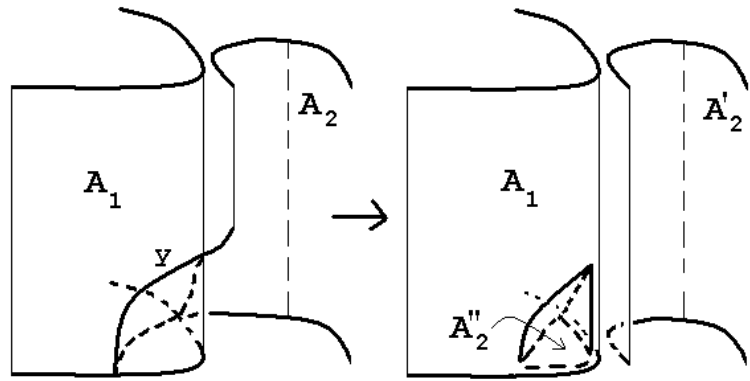}
\caption{Removing the intersection along a non-vertical arc}
\label{rm2}
\end{figure}

Now, let us state two auxiliary lemmas.

\begin{lm}
Suppose $S_{g}$ is the oriented surface of genus $g$ and let
$\Delta$ be  an embedded disk in $S_{g}$. Then if a closed non
self-intersecting curve $\gamma\in S_{g}\backslash \Delta$ is
trivial in $S_{g}$ and not trivial in $S_{g}\backslash \Delta$
then it is parallel to $\partial \Delta$ (i.e. $\gamma\cup
\partial \Delta$ bounds a cylinder in $S_{g}$).\label{lm1}
\end{lm}

Indeed, if a curve $\gamma$  bounds a disk in $S_{g}$; which does
not intersect $\Delta$ then $\gamma$ is contractible in
$S_{g}\backslash \Delta$.

The following lemma is evident.

\begin{lm}
If a proper annulus $C'$ is free homotopic to the annulus $D$
($rel$ boundary),  then $C'$ intersects the knot $K_{1}{\#}K_{2}$.
\label{lm2}
\end{lm}

Now, we may assume that our intersection $C\cap D$ consists only
of vertical arcs, or only of horizontal circles (the existence of
a vertical arc contradicts the existence of a horizontal circle).

Suppose we have only horizontal circles. Then $C$ is homotopic to
$D$ and, by Lemma \ref{lm2}, the annulus $C$ intersects the knot
$K_{1}{\#}K_{2}$. Thus we obtain a contradiction.

Now, suppose the intersection $C\cap D$ consists only of vertical
arcs. Then the annulus $C$ is divided into $2k$ parts
$C_{1},\dots, C_{2k}$, whereas $C_{2l+1}$ lies in $M_{1}\times I$,
and $C_{2l}$ lies in $M_{2}\times I$, where $l\in \{1,\dots, k\}$
when meaningfull. The annulus $C\cap D$ is thus divided into $2k$
sectors by radii (more precisely, radial segments); some of these
sectors contain intersection(s) with the knot, see Fig.
\ref{sect}. Denote all these radii by $r_{1},\dots, r_{2k}$.

\begin{figure}
\centering\includegraphics[width=200pt]{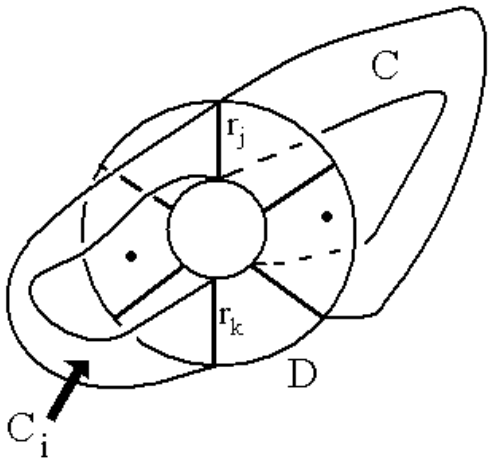} \caption{The
annulus $C$ divided into sectors} \label{sect}
\end{figure}

Now, each part $C_{i}$ of the annulus $C$ is incident to two radii
$r_{j}$ and $r_{k}$. Then $D$ is divided into  two parts
 by these radii. Denote these parts by $D_{jk}^{+}$ and $D_{jk}^{-}$. There are four
options with respect to the following questions:

\begin{enumerate}
\item Is it true that any of the two parts $D_{jk}^{+}$ and
$D_{jk}^{-}$ intersect the knot precisely at one point?

\item Is it true that an annulus obtained by attaching $C_{i}$ to
one of $D_{jk}^{+}$ or $D_{jk}^{-}$, cuts off a ball (so that if
we attach  $C_{i}$ to the other fragment, we get an annulus
homotopic to $D$)?
\end{enumerate}

\begin{re}
Here we mean that a proper surface (say, annulus) $F\subset M$
cuts a ball if $M\backslash F$ has two connected components, one
of which is a topological ball. In other words, $F$ bounds a ball
together with a part $P$ of boundary $\partial M$ of the manifold
$M$, so that $P\cup F$ is a $2$-sphere.
\end{re}

First, consider the case when the answer to the first question is
negative, i.e., one of the parts, say, $D_{jk}^{+}$ does not
intersect the knot (whence the other part $D_{jk}^{-}$ meets the
knot twice).

If the annulus obtained by gluing $C_{i}$ with $D_{jk}^{+}$, cuts
a ball then we may ``pull'' $D_{jk}^{+}$ through $C_{i}$; this
would decrease the number of intersection components between $C$
and $D$. This leads to a contradiction.

If the annulus obtained by pasting $C_{i}$ and $D_{jk}^{-}$ cuts
off a ball then the annulus $C_{i}\cup D_{jk}^{+}$ is homotopic to
 $D$; thus it should intersect the knot. This leads to a contradiction again.

If both answers are affirmative, we get a contradiction: our knot
can not meet the boundary of a ball precisely at one point.

Finally, if, say, $D_{jk}^{+}$ contains precisely one intersection
point with the knot $K_{1}{\#}K_{2}$ then the number of
intersection components of $C\cap D$ should be equal to two.

Indeed, if it is greater than two, then it can be decreased as
shown in Fig. \ref{cnct}. More precisely, among all parts $C_{i}$
of the annulus $C$ we take only two parts and compose a nontrivial
annulus  $C'$, along which $(M_{1}\# M_{2})\times I$ with the knot
$K_{1}{\#}K_{2}$ inside can be destabilized.

In other words, having more than two components $C_{i}$, one can
find two of them which can be repasted and thus obtain a new
non-trivial annulus $C'$ intersecting $D$ at a smaller number of
curves and not intersecting the knot $K_{1}{\#}K_{2}$.

\begin{figure}
\centering\includegraphics[width=400pt]{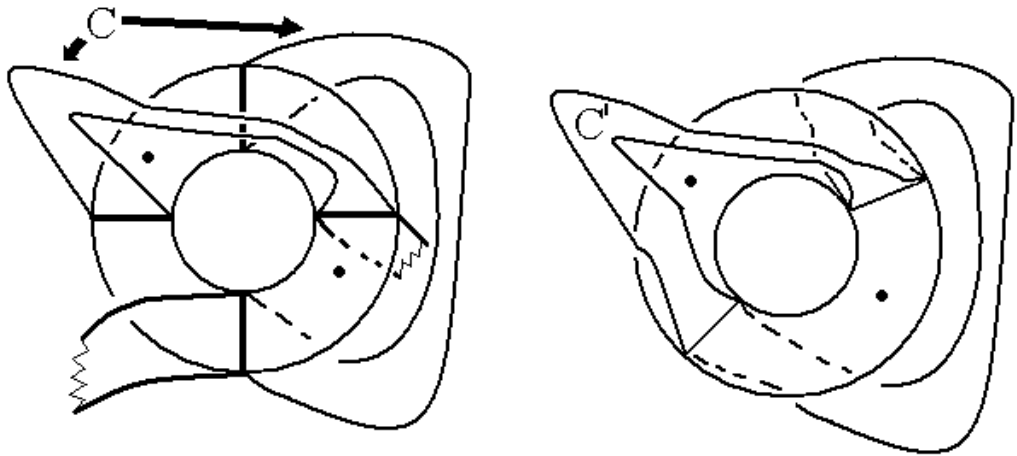}
\caption{Simplifying the curve $C$} \label{cnct}
\end{figure}

Thus we see that if the intersection $C\cap D$ were minimal, there
would be precisely two connected components.

Let us show that the destabilization along such $C$ just
transforms the type of connected sum: we obtain a connected sum of
the second type.

Indeed, the disk $D$ is just cut into two parts by this
destabilization; and two intersection point after the
destabilization lie in different connected components.

The proof for the case when we have the connected sum of the
second type goes in the same vein. Either it is possible to
destabilize only one of $M_{1}$ or $M_{2}$ or all annuli along
which the stabilization can be performed intersect $M_{1}\cap
M_{2}$ (which consists of two components in this case).
Considering such an annulus representing minimal intersection with
$M_{1}\cap M_{2}$, we get the only possibility when the
destabilization transforms the connection type to the first type.
\end{proof}

\section{Algorithmic recognition of virtual links}

The aim of this section is to prove the following

\begin{thm}
There is an algorithm to decide whether two virtual links are
equivalent or not.\label{rcgthm}
\end{thm}

This theorem was first proved in \cite{MProof}, see also
\cite{MBook}.

We shall use the result by Moise \cite{Moi} that each 3-manifold
admits a triangulation. In the sequel, each 3-manifold is thought
to be triangulated.

A manifold $M$ is {\em irreducible} if each embedded sphere in $M$
bounds a ball in $M$.

We shall use the definition of virtual knots as knots in thickened
surfaces ${\cal M}\times I$ up to stabilizations/destabilizations.
Here  ${\cal M}$ is a compact 2-surface, not necessarily
connected. Herewith we require that for each connected component
${\cal M}_{i}$, the 3-manifold ${\cal M}_{i}\times I$ contains at
least one component of the link~$L$.

Recall that a representative for a virtual link is {\em minimal}
if it can not be destabilized.

We shall need Theorem \ref{Kpbg} by Kuperberg \cite{KUP}.

Thus, in order to compare virtual links, we should be able to find
their minimal representatives and compare them. The algorithm to
be given below uses a recognition techniques for three-manifolds
with boundary pattern (see definition below) connected to the
virtual links in question.

We shall use the following facts from Haken-Matveev theory, see
\cite{Mat1}.

A {\em compressing disk} for a surface $F$ in a $3$-manifold $M$
is an embedded disk $ D\subset M$ which meets $F$ along the
boundary of the disk, i.e. $D\cap F =\partial D$.

A surface (possibly disconnected) $F\subset M$ is called {\em
compressible} in one of the two cases:

\begin{enumerate}

\item It admits a compressing disk $D$ such that $\partial D$ does
not bound a disk in $F$;

\item There is a ball in $B$ such that $B\cap F=\partial B$.

\end{enumerate}

A surface is {\em incompressible} if it is not compressible.

A surface $F\subset M$ is called {\em boundary compressible} if
there exists a disk $D^{2}\subset M$ such that $D^{2}\cap
(\partial M\cup F)=\partial D^{2}$ and $D^{2}\cap F$ is a
non-trivial arc in $F$ (an arc that does not cut a disk from $F$).

Also, a 3-manifold $M$ is {\em boundary irreducible} if for any
proper disk $D\subset M$, $\partial D$ bounds a disk on $\partial
M$.

Given a 3-manifold with boundary. By a {\em boundary pattern}
(first proposed by Johannson, see \cite{Joh}) we mean a fixed
$1$--polyhedron (graph) without isolated points on the boundary of
the three manifold (we assume this graph be a subpolyhedron of the
selected triangulation).

The existence of a boundary pattern does not change the definition
of incompressible surface and irreducible manifold.

We have straightforward generalizations of boundary incompressible
surface and boundary irreducible manifold as described below. A
disk $D\subset M$ is called {\em clean} if it does not intersect
the pattern.

For boundary irreducibility we require that every clean proper
disk cuts a  ball (not intersecting the pattern). In the
definition of boundary incompressible surface we require that
every clean disk $D^{2}$ for which $D^{2}\cap (\partial M\cup
F)=\partial D^{2}$ and $l=\partial D^{2}\cap F$ is an arc in $F$,
the arc $l$ cuts a clean disk from $F$.

Let $F$ be a surface in a manifold $(M,\Gamma)$ with boundary
pattern.

Recall that an orientable 3-manifold $M$ is {\em sufficiently
large} if it contains a proper incompressible boundary
incompressible surface distinct from $S^{2}$ and $D^{2}$.

It is natural to consider the notion of sufficiently large
manifold together with properties of irreducibility and boundary
irreducibility. This leads to the following definition.

A connected $3$-manifold without boundary (thus, without boundary
pattern) is {\em Haken} if it is irreducible and sufficiently
large. An irreducible boundary irreducible $3$-manifold
$(M,\Gamma)$ with a boundary pattern $\Gamma$ is Haken either if
it is sufficiently large or if its pattern $\Gamma$ is nonempty
(hence, so is $\partial M$), and  $M$ is a handlebody but not a
ball.

\begin{dfn}
Let $M$ be an irreducible boundary irreducible $3$-manifold. A
proper annulus $A\subset M$ is called {\em inessential} if either
it is parallel $rel$ $\partial$ to an annulus in $\partial M$, or
the core circle of $A$ is contractible in $M$. Otherwise $A$ is
called essential.
\end{dfn}

Any essential annulus in $S_{g}\times I$ with two boundary
components lying in $S_{g}\times \{0\}$ and $S_{g}\times \{1\}$ is
precisely an annulus along which we may destabilize.

A manifold with more than one connected component is called {\em
Haken} if any connected component of it is Haken.

We shall use the following
\begin{prop}[Jaco-Rubinstein-Thompson, see, e.g.,\cite{Tho,Mat1}]
Any connected irreducible  3-manifold with nonempty boundary is
either sufficiently large or a handle body.\label{prp}
\end{prop}

Later on, we deal with manifolds with non-empty boundary pattern.
For this manifold to be Haken, it is sufficient to check that the
manifold in question is irreducible and boundary irreducible but
not a ball.

\begin{lm}[Jaco-Rubinstein-Thompson, see, e.g., \cite{Tho,Mat1}]
There exists an algorithm to decide whether a manifold  $M$ is
reducible; if it is so, the algorithm constructs a $2$--sphere
$S\subset M$ not bounding a ball in $M$.\label{pRiv}
\end{lm}

\begin{lm}\(\cite{Mat1}\)
Classical links are algorithmically recongizable.
\end{lm}

This lemma follows from Haken's theory of normal surfaces; the
proof is based on the following ideas: for each non-trivial
non-split link, the complement in $S^{3}$ to the tubular
neighborhood of this link is a Haken manifold. Endowing the
boundary with a pattern, we will be able to restore the initial
link. After that, the problem is reduced to the recognition
problem for Haken manifolds, for details see \cite{Mat1}.

\begin{lm}\(\cite{Mat1}\)
There is an algorithm to decide whether a Haken manifold $M$ with
a boundary pattern $\Gamma\subset \partial M$ has a proper clean
essential annulus. If such an annulus exists, it can be
constructed algorithmically.\label{suKol}
\end{lm}

We shall use this lemma to define whether a given representative
of a virtual link can be destabilized.

\begin{lm}\(\cite{Mat1}\)
There is an algorithm to decide, whether two Haken manifolds
$(M,\Gamma)$ and $(M',\Gamma')$ with boundary patterns are
homeomorphic by means of a homeomorphism that maps $\Gamma$ to
$\Gamma'$.\label{Hra}
\end{lm}

Consider an arbitrary representative of a virtual link $L$, i.e.,
a couple $(M,L)$, where $M={\tilde M}\times I$ for some closed
2-surface ${\tilde M}$, and  $L$ is a link in $M$ (we use the same
letter  $L$ for denoting both the initial link and the
representing link in $M$: abusing notation). Let $N$ be a small
open tubular neighborhood of the link $L$. Cut  $N$ from $M$. We
obtain a manifold with boundary. Denote it by $M_{L}$. Its
boundary consists of boundary components of $M$ (two, if  $M$ is
connected) and several tori; the number of tori equals the number
of components of the link $L$. Let us endow each torus with a
pattern $\Gamma_{L}$, representing the meridian of the
corresponding component (we also add a vertex to make a graph from
the meridional circle). Thus we obtain the manifold
$(M_{L},\Gamma_{L})$ with a boundary pattern.

It is obvious that the virtual link  $L$ (and the pair $(M,L)$)
can be restored from $(M_{L},\Gamma_{L})$, since we know how to
restore the manifold $M$ by attaching full tori to the boundary
components of  $M_{L}$ knowing meridians of these full tori.

\begin{lm}
Suppose a link $L$ is not a split sum of a (nonempty) classical
link and a virtual link. Then the manifold $(M_{L},\Gamma_{L})$
with boundary pattern $\Gamma_{L}$ is Haken.\label{V6}
\end{lm}

\begin{proof}
In virtue of Proposition \ref{prp}, it remains to show that (any
connected component) this manifold with boundary pattern is
irreducible and boundary irreducible: by definition, it can not be
a handlebody.

In the case when $g=0$ we deal with classical links. Suppose
$g>0$. Then, for any connected orientable $2$-surface $S_{g}$ the
manifold $S_{g}\times I$ is irreducible.  Thus, if the link  $L$
is not classical (thus $g\neq 0$) then for its neighborhood
$N(L)$, the set $(S_{g}\times I)\backslash N(L)$ might be
reducible if and only if it contains a sphere  $S$, bounding a
ball in $S_{g}\times \{0,1\}$ such that this ball contains some
components of the link $L$. This means that these components form
a classical sublink of $L$ separated from all other components.

Furthermore, since $L$ is not a split sum of the unknot with some
virtual link, the manifold $M_{L}$ is boundary irreducible.

Indeed, each curve in   $S_{g}\times \{0\}$ or in $S_{g}\times
\{1\}$ which may bound a disk in $S_{g}\times I$ is contractible
in the boundary. Thus, boundary reducibility can occur only if we
have a proper disk with boundary lying on some torus --- the
boundary of the cut full torus. This should mean that the cut full
torus corresponds to the split unknot of the link.

Thus, the considered manifold is irreducible and boundary
irreducible and thus (by Proposition \ref{prp}), Haken.
\end{proof}

Now, let us prove the main theorem. Let  $L,L'$ be virtual links.

\begin{enumerate}

\item[Step 1.] Consider some representatives $(M,L), (M',L')$ of
the virtual links in question. Let us construct the corresponding
manifolds with boundary patterns. Denote them by
$(M_{L},\Gamma),(M'_{L'},\Gamma')$.

\item[Step 2.] Define whether one of $M_{L}$ or $M'_{L'}$ is
reducible. If one of them is so, then, by Lemma  \ref{pRiv}, we
may find a sphere not bounding a ball, and thus separate some
classical components of the corresponding link.

Now, rename the manifolds with boundary patterns accordingly: i.e.
we shall use the previous notation for what is left from manifolds
by cutting off classical components.

\item[Step 3.] Define (by Lemma \ref{suKol}) whether it is
possible to destabilize one of  $(M,L)$ or $(M',L')$. If it is
possible, perform the destabilization. Return to Step 2.

Let us perform steps 2 and 3 while possible. Obviously, this
process stops in a finite period of time.

Classical links are algorithmically recognizable. Thus, we may
compare the split classical sublinks of  $L$ and $L'$. If they are
not isotopic, we stop: the virtual links in question are not
equivalent. Otherwise, we go on.

After performing the first three steps, we reduce our problem to
the case when there are no split components and representatives
are minimal. From now on, the manifolds in question are Haken by
Lemma \ref{V6}.

\item[Step 4.]  Each connected component of the manifolds
$(M_{L},\Gamma)$ and $(M'_{L'},\Gamma')$ is a Haken manifold with
a boundary pattern. Thus, we can algorithmically solve the problem
whether there exists a homeomorphism $f:M_{L}\to M'_{L'}$ that
maps $\Gamma$ to $\Gamma'$ (by Lemma \ref{Hra}). If such a
homeomorphism exists then virtual links $L,L'$ are equivalent.
Otherwise $L$ and $L'$ are not equivalent.

\end{enumerate}

Performing the steps described above, we solve the recognition
problem. Theorem \ref{rcgthm} is proved.

\begin{re}
The proof given above works also for oriented virtual links and
framed virtual links.
\end{re}

\section{Self-linking Numbers for Virtual Links}
Call a classical  crossing  in an oriented virtual knot diagram
$K$ {\em odd} if, in the Gauss code for that diagram there are an
odd number of appearances of (classical) crossings between the
first and the second appearance of $i.$ Let
$$J(K) = w(K)|_{Odd(K)}$$ where $Odd(K)$ denotes the collection of odd crossings of $K,$ and the restriction of the writhe
to $Odd(K)$, $w(K)|_{Odd(K)},$ means the summation over the signs of the odd crossings in $K.$ Then
it is not hard to see that $J(K)$ is invariant of the virtual knot or link $K.$ We call $J(K)$ the {\em self-linking number}
of the virtual diagram $K.$ This invariant is simple, but remarkably powerful.
\bigbreak

\noindent

If $K$ is classical then $J(K) = 0,$ since there are no odd
crossings in a classical diagram.

\begin{thm} Let $K$ be a virtual knot and let $K^{*}$ denote the
mirror image of $K$ (obtained by switching all the crossings of
the diagram $K$). Then
$$J(K^{*}) = -J(K).$$\label{th4}
\end{thm}

Hence, if $J(K)$ is non-zero, then $K$ is inequivalent to its
mirror image. If $K$ is a virtual knot and $J(K)$ is non-zero,
then $K$ is not equivalent to a classical knot.

\bigbreak

We leave the proof of this Theorem and the proof of the invariance of $J(K)$ to the reader.
See \cite{SelfLink} for more about this invariant its generalizations.
\bigbreak

View Figure \ref{S}. The two virtual knots in this figure
illustrate the application of Theorem \ref{th4}. In the case of
the virtual trefoil $K,$ the Gauss code of the shadow of $K$ is
$abab;$ hence both crossings are odd, and we have $J(K) = 2.$ This
proves that $K$ is non-trivial, non-classical and inequivalent to
its mirror image. Similarly, the virtual knot $E$ has shadow code
$abcbac$ so that the crossings $a$ and $b$ are odd. Hence $J(E) =
2$ and $E$ is also non-trivial, non-classical and chiral. Note
that for $E,$ the invariant is independent of the type of the even
crossing $c.$ \bigbreak

\begin{figure}
\centering\includegraphics[width=200pt]{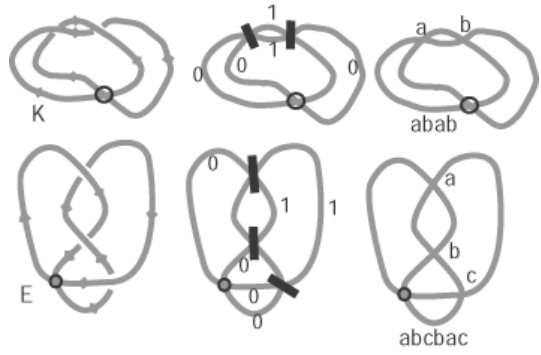}
\caption{Virtual Trefoil $K$ and Virtual Figure Eight $E$}
\label{S}
\end{figure}

View Figure \ref{abcda}. The virtual knot $K'$ in this figure has
Gauss code $abcdcadb$, and $J(K')=2$. Note that $K'$ would be
unknotted if we allowed the forbidden move (of other type). This
example underlines why we forbid such moves in virtual knot
theory.

\begin{figure}
\centering\includegraphics[width=200pt]{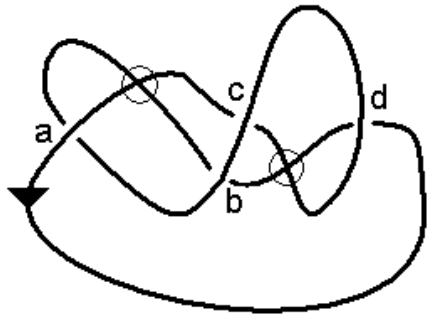} \caption{The
Knot $K'$} \label{abcda}
\end{figure}

\section{Welded Braids and Tubes in Four-Space}

The welded braid group $WB_{n}$ can be interpreted as the fundamental group of the configuration space of
 $n$ disjoint circles trivially embedded in three dimensional space ${\bf R}^{3}$. This group (the so-called motion
group of disjoint circles) can, in turn, be interpreted as a braid
group of tubes embedded in ${\bf R}^{3} \times {\bf R} = {\bf
R}^{4}.$ These braided tubes in four-space are generated by two
types of elementary braiding. In Figure \ref{21}, we show diagrams
that can be interpreted as immersions of tubes in three-space.
Each such immersion is a projection of a corresponding embedding
in four-space. The first two diagrams of Figure \ref{21} each
illustrate a tube passing through another tube. When tube $A$
passes through tube $B$ we make a corresponding classical braiding
crossing with arc $A$ passing under arc $B.$ The four-dimensional
interpretation of tube $A$ passing through tube $B$ is that:
 As one looks at the levels of intersection with
${\bf R}^{3} \times t$ for different values of $t,$ one sees two
circles $A(t)$ and $B(t).$ As the variable $t$ increases, the
$A(t)$ circle (always disjointly embedded from the $B(t)$ circle)
moves through the $B(t)$ circle. This process is illustrated in
Figure \ref{22}. \bigbreak

While the classical crossing in a welded braid diagram corresponds
to a genuine braiding of the tubes in four-space (as described
above), the virtual crossing corresponds to tubes that do not
interact in the immersion representation  (see again Figure
\ref{21}). These non-interacting tubes can pass over or under each
other, as these local projections correspond to equivalent
embeddings in four-space. \bigbreak

\begin{figure}
\centering\includegraphics[width=200pt]{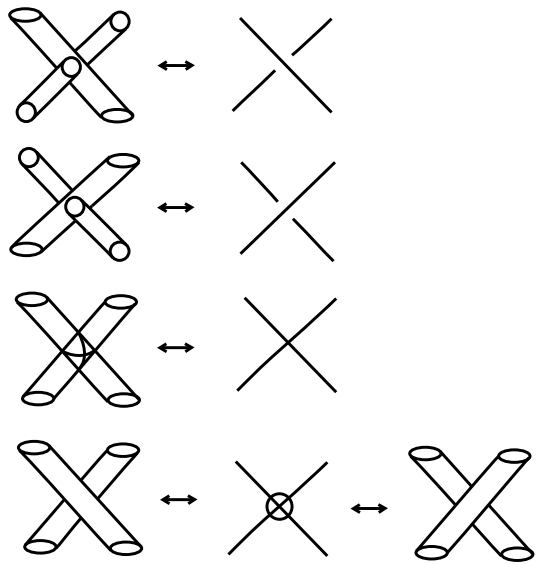}
\caption{Tubular
correspondence} \label{21}
\end{figure}

\begin{figure}
\centering\includegraphics[width=200pt]{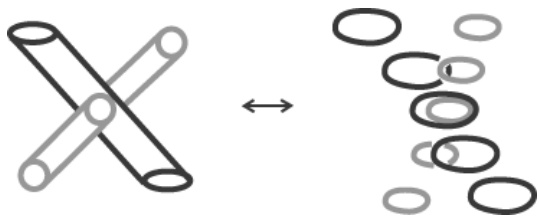}
\caption{Braiding of circles} \label{22}
\end{figure}

It is an interesting exercise to verify that the moves in the
welded braid group each induce equivalences of the corresponding
tubular braids in four-space. In particular, the move $(F_1)$
induces such an isotopy, while the forbidden move $(F_2)$ does
not. For more on this subject, the reader can consult \cite{Satoh}
and also \cite{KaV2} and the references therein. The basic idea
for this correspondence is due to Satoh in \cite{Satoh} where
torus embeddings in four-space are associated with virtual knot
diagrams. \bigbreak

For knot theory the moral of these remarks is that the category of
welded knots and links (virtual knots and links plus the
equivalence relation generated by the first forbidden move) is
naturally associated with embeddings of tori in four dimensional
space. To each welded knot or link there is associated a
well-defined embedding of a collection of tori (one torus for each
component of the given link) and the fundamental group of the
complement of this embedding in four space is isomorphic to the
combinatorial fundamental group of the welded link (this is the
same as the combinatorial fundamental group of a corresponding
virtual link). It is an open problem whether this association is
an embedding of the category of welded links into the category of
toroidal embeddings in four-space. \bigbreak

\end{document}